\documentclass[final,1p,times]{elsarticle}

\usepackage{amssymb}

\usepackage{amsmath}

\usepackage{caption}
\usepackage{subfig}
\usepackage{bm}
\usepackage{booktabs}
\usepackage[numbers]{natbib}
\usepackage{url}  
\usepackage{float}
\usepackage{hyperref}
\usepackage{color}
\usepackage[section]{placeins}
\journal{}

\begin{document}
	
	\begin{frontmatter}
		
		\title{Data-driven optimized high-order WENO schemes with low-dissipation and low-dispersion}
		
		\author[a]{Jinrui Zhou}
		\ead{zhou1218@std.uestc.edu.cn}
		
		\author[a]{Yiqi Gu}
		\ead{yiqigu@uestc.edu.cn}
		
		\author[b]{Song Jiang}
		\ead{jiang@iapcm.ac.cn}
		
		\author[a]{Hua Shen}
		\ead{huashen@uestc.edu.cn}
		
		\author[a]{Liwei Xu\corref{cor1}}
		\ead{xul@uestc.edu.cn}
		
		\author[c]{Guanyu Zhou}
		\ead{zhoug@uestc.edu.cn}

		\affiliation[a]{
			organization={School of Mathematical Sciences, University of Electronic Science and Technology of China},
			city={Chengdu},
			state={Sichuan},
			postcode={611731},
			country={China}
		}
		\affiliation[b]{
			organization={National Key Laboratory of Computational Physics, Institute of Applied Physics and Computational Mathematics},
			city={Beijing},
			postcode={100094},
			country={China}
		}
		\affiliation[c]{
			organization={Institute of Fundamental and Frontier Sciences, University of Electronic Science and Technology of China},
			city={Chengdu},
			state={Sichuan},
			postcode={611731},
			country={China}
		}
		\cortext[cor1]{Corresponding author}

	\begin{abstract}
		Classical high-order weighted essentially non-oscillatory (WENO) schemes are designed to achieve optimal convergence order for smooth solutions and to maintain non-oscillatory behaviors for discontinuities. However, their spectral properties are not optimal, which limits the ability to capture high-frequency waves and small-scale features.
		In this paper, we propose a data-driven optimized method to improve the spectral properties of the WENO schemes. 
		By analyzing the approximate dispersion relation (ADR), the spectral error of the schemes can be bounded by the reconstructed errors of a series of trigonometric functions with different wavenumbers. 
		Therefore, we propose the new schemes WENO5-JS/Z-NN that introduce a compensation term parameterized by a neural network to the weight function of the WENO5-JS/Z schemes. 
		The neural network is trained such that the generated weights can minimize the reconstructed errors over a large number of spatial stencils, and furthermore, improve the spectral accuracy. 
		Meanwhile, the Total Variation Diminishing (TVD) constraint and anti-dissipation penalization are incorporated into the loss function to enhance the shock-capturing capability and preserve stability in simulating high-frequency waves.
		Compared to WENO5-JS/Z, our schemes maintain the ability to capture discontinuities while providing higher resolution for fine-scale flow features.
		The ADR indicates that the new schemes can match the exact spectrum more accurately over a broader range of wavenumbers.
	\end{abstract}
	
	
	\end{frontmatter}
	
	\section{Introduction}
	\label{section1}
	High-fidelity simulation for hyperbolic conservation laws remains challenging due to the coexistence of discontinuities and multiscale structures in solutions. 
	These intrinsic features impose conflicting demands on numerical schemes: achieving high resolution requires low numerical dissipation and dispersion, whereas suppressing spurious oscillations near discontinuities necessitates sufficient dissipation to maintain stability. 
	In the past decades, various numerical methods have been proposed, including monotonicity-preserving schemes \cite{ref1, ref2}, TVD schemes \cite{ref3, ref4}, essentially non-oscillatory (ENO) schemes \cite{ref5, ref6, ref7}, and WENO schemes \cite{ref7, ref8}, which aim to achieve both high-order accuracy and robust shock-capturing capability.
	
	Liu et al. \cite{ref8} first proposed WENO schemes by using a convex combination of local candidate stencils that exist in ENO schemes. 
	The original WENO schemes are improved by Jiang and Shu \cite{ref9} through introducing a new approach to compute smoothness indicators and nonlinear weights for the candidate stencils, enabling the schemes to achieve optimal convergence order in smooth regions.
	Henrick et al. \cite{ref10} identified that WENO schemes proposed by Jiang and Shu lose accuracy at critical points and addressed this issue by applying a mapping function to the nonlinear weights.
	Borges et al. \cite{ref11} designed a higher-order global smoothness indicator and nonlinear weight function to handle this question. Up to now, the extended WENO schemes, such as the monotonicity-preserving WENO schemes \cite{ref12}, the weighted compact nonlinear schemes \cite{ref13, ref14}, and the central WENO schemes \cite{ref15, ref16} have been developed.   
	
	To explore low-dissipation and low-dispersion WENO schemes for turbulence and aeroacoustic simulations, extensive research has focused on the background linear scheme and the nonlinear weight mechanism. 
	The WENO-SYMBO \cite{ref17} added a downwind stencil to the candidates and optimized the background linear scheme for better spectral properties. 
	Sun et al. \cite{ref19} derived a sufficient condition for the linear scheme to have independent dispersion and dissipation, and designed a class of schemes with minimized dispersion and controllable dissipation (MDCD), where the dissipation parameter is adjusted based on experience in practice. 
	The MDCD is then extended by Hu et al. \cite{ref20} to a sixth-order scheme using the dispersion-dissipation condition \cite{ref21} to determine an adequate dissipation parameter. Moreover, Li et al. \cite{ref22} developed a method that can adaptively adjust the dissipation parameter according to the local flow scale. 
	Recently, Acker et al. \cite{ref24} and Luo et al. \cite{ref25, ref26} improved the numerical resolution of the smooth waves by increasing the contribution of the sub-stencils where the solution is less smooth. The spectral properties of these new methods are superior to the linear scheme. 
	A family of high-order targeted ENO (TENO) schemes was constructed by Fu et al. \cite{ref27, ref28}, where the ENO-like stencil selection strategy significantly reduced numerical dissipation.
	
	In recent years, data-driven methods based on fully connected neural networks (FCNNs) have gained significant attention in computational fluid dynamics due to their ability to model complex nonlinear relationships present in data. 
	They have been integrated with classical numerical schemes to improve numerical accuracy or computational efficiency, including accelerating high-order DG solvers \cite{ref29, ref30}, surrogating the smoothness detector \cite{ref31}, and constructing high-precision schemes on coarse grids  \cite{ref32, ref33, ref34}. 
	Stevens and Colonius \cite{ref35} pioneered the incorporation of machine learning into WENO schemes to enhance numerical resolution, while the accuracy is reduced to first order. 
	To address this limitation, the WENO3-NN method \cite{ref36} was proposed, which introduces the adaptive penalization to the loss function, and successfully restores the maximum-order convergence in smooth regions. Nogueira et al. \cite{ref37} further extended this strategy to the fifth-order scheme. 
	Fan et al. \cite{ref41} introduced a multi-resolution strategy to mitigate the issue of accuracy degradation observed in the WENO3-NN method as the mesh refines. 
	A general deep-reinforcement-learning framework is developed by Feng et al. \cite{ref38} to find the optimal parameters used in the fifth-order TENO. 
	Liu et al. \cite{ref39} trained a neural network to adaptively regulate the dispersion and dissipation coefficients of the six-point linear schemes. 
	A spectral optimization framework based on neural networks for the numerical schemes on unstructured grids is introduced by Zhou et al. \cite{ref40}. 
	For more recent advances in neural networks for computational fluid dynamics, one may refer to the review \cite{ref42}.
	
	In this work, we provide a data-driven optimized method WENO5-NN to improve the spectral properties of the classical fifth-order WENO-JS/Z schemes.
	The strategy is to introduce a compensation term to the nonlinear weights function of the original schemes, which is parameterized by a neural network and learned by the data-driven method. 
	Through the analysis of the ADR, we find that for any fixed scaled wavenumber, the spectral error can be bound by the reconstructed errors of the trigonometric functions at that wavenumber.
	By minimizing the reconstruction error between numerical and reference interface fluxes over a broad range of spatial stencils containing sine functions, the spectral accuracy of the WENO5-JS/Z-NN schemes can be enhanced.
	The TVD condition and the anti-dissipation penalization are incorporated into the loss function to enhance the ability to capture discontinuities and maintain stability in high wavenumbers, respectively. 
	A series of experiments is used to demonstrate the validity of the penalization terms. 
	Compared with the nonlinear weights of WENO5-JS/Z, the new schemes adaptively adjust the background linear weights for the smooth wave and WENO5-Z-NN slightly increases the contribution of the less-smooth stencils near discontinuities. 
	The ADR reveals that the proposed schemes have lower dispersion and dissipation than WENO5-JS/Z in the mid-frequency wavenumbers.
	Several benchmark cases of the Euler equations demonstrate that the WENO5-JS/Z-NN schemes preserve robust shock-capturing capabilities while achieving high resolution.

	The rest of this paper is organized as follows. Section \ref{section2} reviews the classical WENO5-JS and WENO5-Z schemes. Section \ref{section3} introduces the spectral error estimate and the optimized WENO method. The comparison of the numerical examples between the WENO5-JS/Z and WENO5-JS/Z-NN schemes is presented in Section \ref{section4}. Finally, we give a conclusion in Section \ref{section5}.
	
	\section{Preliminary}
	\label{section2}
	In this section, we review the basic concepts of classical WENO schemes. The reconstructed process of the fifth-order WENO-JS and WENO-Z is then presented.
	
	\subsection{Basic concepts of WENO schemes}
	\label{section2.1}
\begin{figure}
		\centering
		\includegraphics[scale=0.45]{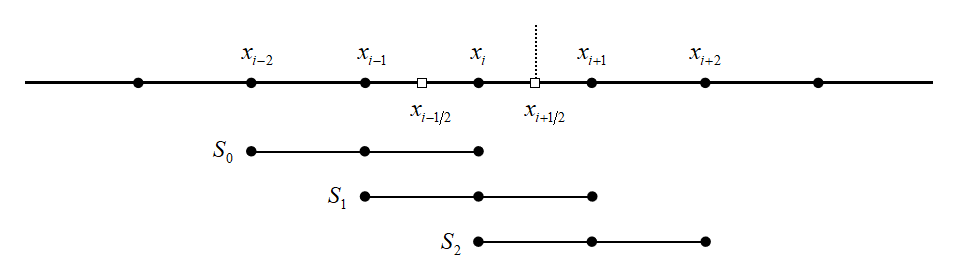}
		\caption{Candidate stencils for the WENO-JS and WENO-Z schemes.}
		\label{fig:2.1}
	\end{figure}

	Here, we consider the one-dimensional scalar conservation law
	\begin{equation}
		\frac{\partial u}{\partial t} + \frac{\partial  f(u)}{\partial x} = 0, x\in[x_L, x_R], t \in[0, \infty),
		\label{eq:2.0.1}
	\end{equation}
	where $u$ is the solution and $f$ is the flux function. We assume $\frac{\partial f}{\partial x} \ge 0$ without loss of generality. The spatial domain $[x_L, x_R]$ is uniformly discretized into some interval $[x_{i-\frac{1}{2}}, x_{i+\frac{1}{2}}]$, $i = 1, \dots, N_x$, where $N_x$ is the number of mesh grids and the cell center is $x_i = \frac{1}{2}(x_{i-\frac{1}{2}} + x_{i+\frac{1}{2}})$. $u_i(t)$ is the numerical approximation to the nodal point value $u(x_i, t)$. Eq. (\ref{eq:2.0.1}) can be written as a conservative finite difference scheme
	\begin{equation}
		\frac{d u_i(t)}{d t} = -\frac{h_{i+\frac{1}{2}} -  h_{i-\frac{1}{2}}}{\Delta x},
		\label{eq:2.0.2}
	\end{equation} 
	where $h(x)$ is implicitly defined as $f(x) = \frac{1}{\Delta x} \int_{x - \frac{\Delta x}{2}}^{x + \frac{\Delta x}{2}} h(\xi) ~d\xi$ and $h_{i\pm\frac{1}{2}} = h(x_{i\pm\frac{1}{2}})$.
	Eq. (\ref{eq:2.0.2}) can be approximated as
	\begin{equation}
		\frac{d u_i(t)}{d t} \approx  -\frac{\hat{f}_{i+\frac{1}{2}} - \hat{f}_{i-\frac{1}{2}}}{\Delta x},
		\label{eq:2.0.3}
	\end{equation}
	where the numerical flux $\hat{f}_{i \pm \frac{1}{2}}$ is high-order approximation to the interface flux $h_{i\pm\frac{1}{2}}$. The fifth-order WENO flux is reconstructed by a convex combination of the three low-order numerical fluxes on the substencils $S_0$, $S_1$, $S_2$, where the substencils are shown in Fig. \ref{fig:2.1}. The low-order numerical fluxes are defined as
	\begin{align*}
			&\hat{f}_{i+\frac{1}{2}}^0 = \frac{1}{3}  f_{i-2} - \frac{7}{6}  f_{i - 1} + \frac{11}{6} f_{i}, ~
			\hat{f}_{i+\frac{1}{2}}^1 =  -\frac{1}{6}f_{i-1} +  \frac{5}{6}f_{i} + \frac{1}{3} f_{i+1},  \\
			&\hat{f}_{i+\frac{1}{2}}^2 = 	 \frac{1}{3}f_{i} +  \frac{5}{6}f_{i+1} - \frac{1}{6} f_{i+2}. \quad
	\end{align*}
	The final flux $\hat{f}_{i+1/2}$ is then given by
	\begin{equation}
		\hat{f}_{i+\frac{1}{2}} = \sum_{k = 0}^{2} w_k  \hat{f}_{i+\frac{1}{2}}^k.
		\label{eq:2.0.5}
	\end{equation}
	
	The design of the nonlinear weights $w_k$ should satisfy the following two conditions: 
	\begin{enumerate}
		\item When a discontinuity is present in the global stencil $\{ S_0, S_1, S_2 \}$, the weight associated with the substencil containing the discontinuity should be sufficiently small. 
		\item When the global stencil is smooth, the nonlinear weights should recover higher-order accuracy as much as possible.
	\end{enumerate}
	
	\subsection{Fifth-order WENO-JS and WENO-Z schemes}
	\label{section2.2}
	The WENO-JS weights \cite{ref9} are defined as
	\begin{equation}
		w_k^{JS} =  \frac{\alpha_k}{\sum_{l=0}^{2}\alpha_l}, \quad \alpha_k =  \frac{d_k}{(\beta_k + \epsilon)^p}, \quad k = 0, 1, 2, 
		\label{eq:2.1.1}
	\end{equation} 
	where $(d_0,d_1,d_2)=(0.1, 0.6, 0.3)$ are the ideal weights used to generate a fifth order upwind central scheme (UP5) on the global stencil, $\epsilon$ is a small value to avoid division by zero. The smoothness indicators $\beta_k$ measure the regularity of the linear polynomial at the stencil $S_k$, which is defined as 
	\begin{subequations}
		\begin{align}
			& \beta_0 = \frac{13}{12}(f_{i-2} - 2  f_{i-1} + f_{i})^2	 +  \frac{1}{4}(f_{i-2} - 4 f_{i-1} + 3f_{i})^2 \label{eq:2.1.2a}, \\
			& \beta_1 = \frac{13}{12}(f_{i-1} - 2 f_{i}  + f_{i+1})^2 +   \frac{1}{4}(f_{i-1} - f_{i+1})^2  \label{eq:2.1.2b}, \\
			& \beta_2 = \frac{13}{12}(f_{i} - 2 f_{i+1}  + f_{i+2})^2 +  \frac{1}{4}(3f_{i} - 4f_{i+1} + f_{i+2})^2 \label{eq:2.1.2c}.
		\end{align}
	\end{subequations}	
	The power $p$ controls the sensitivity of the nonlinear weights to smoothness indicators and is usually set to $p = 2$.
	
	In the smooth region, the Taylor expansion for the smoothness indicators shows that $\beta_k = O(\Delta x^2)$, implying that $w_k = d_k + O(\Delta x^2)$. Therefore, the nonlinear weights $w_k$ are close to the ideal weights. When the stencil $S_k$ contains discontinuity, $\beta_k$ is $O(1)$ while other indicators are $O(\Delta x^2)$. The corresponding weight $w_k$ is very small and can avoid numerical oscillation. 
	
	Borges et al. \cite{ref11} investigated the	effective order of classical WENO schemes and provided a sufficient condition for fifth-order convergence. They proposed the WENO-Z scheme:   
	\begin{equation}
		w_k^{Z} =  \frac{\alpha_k}{\sum_{l=0}^{2}\alpha_l},  \quad \alpha_k =  d_k\left(1 + \left(\frac{\tau_5}{\beta_k +  \epsilon}\right)^q\right), \quad k = 0, 1, 2,
		\label{eq:2.1.2}
	\end{equation} 
	where $\tau_5 = |\beta_2 - \beta_0|$ is the global smoothness measure and the power is suggested to be set to $q = 1$. If the global stencil does not contain discontinuity, then $\tau_5 = O(\Delta x^ 5)$ and the nonlinear weights satisfy $w_k^Z = d_k + O(\Delta x^3)$, $k = 0, 1, 2$. Compared to WENO-JS, WENO-Z exhibits lower numerical dissipation as its weights are closer to the optimal values $d_k$.
	
	\section{Optimized WENO schemes with deep learning}
	\label{section3}
	In this section, we first introduce some fundamental concepts, including general neural networks and the TVD condition. We then provide a spectral error estimate for nonlinear schemes, which reveals that one can improve the spectral accuracy by minimizing the reconstruction error of the learned weights. Finally, the detailed construction of the optimized WENO scheme is presented.
	
	\subsection{Neural networks}
	\label{section3.1}
	Neural networks \cite{ref43} represent a class of complex function spaces with parameters. For any input $\bm{x}$, the output value $\bm{y}$ of the neural network with $L$-layers is defined as:
	\begin{equation}
		\begin{aligned}
			& \text{Input layer}  &
			& \bm{x}^{[0]} = \bm{x},	\\
			& \text{Hidden layer} &
			& \bm{x}^{[l]} = \sigma\left( \bm{W}^{[l]}   \bm{x}^{[l-1]} +  \bm{b}^{[l]} \right), \quad l= 1, \dots, L-1,	\\
			& \text{Output layer} & 
			& \mathbf{y} = \sigma\left( \bm{W}^{[L]}   \bm{x}^{[L-1]} +  \bm{b}^{[L]} \right),
			\label{eq:3.1.1}
		\end{aligned}
	\end{equation}
	where $\bm{W}^{[l]}$ and $\bm{b}^{[l]}$ are the weight matrix and bias vector of the neural network at layer $l$, and $\sigma(\cdot)$ is the activation function on every element of the vector. In general, the activation function in the output layer is the identity function. The difference between the output of the neural network $\bm{y}$ and the target value $\hat{\bm{y}}$ is quantified by a loss function, denoted as $\mathcal{L}(\bm{y}, \hat{\bm{y}})$. The training process is to find a set of parameters $\theta$ such that it approximately minimizes the function $\mathcal{L}(\bm{y}, \hat{\bm{y}})$.
	
	\subsection{Total Variation Diminishing condition}
	\label{section3.2}
	The TVD condition is an essential concept in the design of numerical schemes for solving hyperbolic conservation laws. Its primary objective is to suppress non-physical oscillations in numerical solutions, particularly near discontinuities, thereby ensuring the numerical stability of the solution.

	For a discrete solution $u_n := \{u_{i, n}\}_i$, the total variation (TV) is defined as 
	\begin{align}
		\text{TV}(u_n) = \sum_{i}^{}|u_{i+1, n} - u_{i,  n}|.
		\notag
	\end{align}
	TV measures the degree of the total oscillation of the numerical solution in the spatial direction. A scheme is TVD if
	\begin{align}
		\text{TV}(u_{n+1}) \leq \text{TV}(u_{n}) \label{eq:3.2.7a},
	\end{align}
	for all time steps $n$, which indicates that the degree of the total oscillation of the numerical solution does not increase with temporal evolution. In this work, we append the TVD condition to the optimized WENO method to enhance the numerical robustness near discontinuities.
	
	\subsection{Spectral properties of nonlinear schemes}
	\label{section3.3}
	Spectral analysis plays a crucial role in developing low-dispersion and low-dissipation WENO schemes. In this section, we adopt the approximate dispersion relation (ADR) introduced by Pirozzoli \cite{ref44} to evaluate the spectral properties of nonlinear schemes.
	The spectrum obtained in the ADR depends on the time step size $\tau$, suggesting that $\tau$ should be chosen sufficiently small to restrain the influence of temporal discretization. To analyze the effects of spatial discretization errors on the spectrum of nonlinear schemes, we derive a time-independent form of the ADR.

	To calculate the numerical spectrum of the investigated scheme, we first set up a domain with length $L$ and the nodes $x_j = j\delta$, $j = 0, \dots, N$, where $\delta = L/N$. On the discrete grid, the supported Fourier modes have wavelengths $\lambda_n = L / n$, $n = 0, \dots, N/2$. We consider the one-dimensional linear advection equation 
	\begin{align*}
		\frac{\partial u}{\partial t} + \frac{\partial f}{\partial x} = 0, \quad f(u) = au.
	\end{align*}
	For simplicity, we assume $a>0$. A pure harmonic wave $u(x) = \mathcal{F}[u](w_n) e^{iw_nx}$ is used as the initial condition, where $w_n = 2\pi / \lambda_n$ is the discrete wavenumber and $\mathcal{F}[u](w_n)$ is the Fourier coefficient of $u(x)$ at wavenumber $w_n$. The corresponding reduced wavenumber is $\varphi_n = w_n \delta$. The discrete spatial derivative $f^\prime(x_j)$ is exactly expressed as 
	\begin{align*}
		f^\prime (x_j) = \frac{i a \varphi_n \mathcal{F}[u](w_n)}{\delta} e^{i\varphi_n j}: = \mathcal{F}[f^\prime](\varphi_n) e^{i\varphi_n j}.
	\end{align*}
	We can also define a implicit function $h(x)$ to represent $f^\prime (x_j)$ as $f^\prime (x_j) = (h_{j+1/2} - h_{j-1/2})/\delta$, where $ h_{j\pm1/2}:= h(x_j \pm \delta/2)$.

	As discussed in section \ref{section2.1}, the spatial derivative is approximated by classical WENO schemes as 
	\begin{align*}
		f^\prime_j = \frac{1}{\delta}(\hat{f}_{j+1/2} - \hat{f}_{j-1/2}),
	\end{align*}
	where $\hat{f}_{j+1/2} = \hat{f}(f_{j-q+1}, \dots, f_{j+r})$. $\hat{f}$ is a nonlinear function of its argument and therefore no analytical formula of the modified wavenumber can be obtained in this case. 
	By calculating the DFT coefficient $\mathcal{F}[f^\prime_j](\varphi_n)$ of $f^\prime_j$ at $\varphi_n$, the numerical modified wavenumber $\Phi(\varphi_n)$ is obtained as 
	\begin{align}
		\Phi(\varphi_n) = -i\frac{ \mathcal{F}[f^\prime_j](\varphi_n) \delta}{ a \mathcal{F}[u](w_n)}. 
		\label{eq:2.2.1} 
	\end{align}
	The real and imaginary parts of $\Phi(\varphi) - \varphi$ represent the numerical dispersion and dissipation of the nonlinear scheme, respectively. 
	
	We first give the mathematical representation of the core variables to estimate the spectral error $\left|\Phi(\varphi_n) - \varphi_n\right|$. Let $\mathcal{F}[u](w_n) = 1$ for simplicity. Defining the $u(x) = e^{iw_n x} = \cos(w_n x) + i \sin(w_n x): = u^R(x) + i u^I(x)$, $f^\prime(x)$ can be written as $f^\prime_{R}(x) + i f^\prime_{I}(x)$. Let the symbol $*$ stand for either $R$ or $I$.
	By the definition of $h(x)$, we have 
	\begin{align}
		f^\prime_{*}(x_j) = \frac{(h_{*, j+1/2} - h_{*, j-1/2})}{\delta}.
		\label{eq:A.1}
	\end{align}
	The numerical derivative $f^\prime_j$ is decomposition as $f^\prime_{R, j} + i f^\prime_{I, j}$, where $f^\prime_{R, j}$ and $f^\prime_{I, j}$ are the numerical derivative for $f_R(x)$ and $f_I(x)$ at $x = x_j$. By the definition of $\hat{f}_{j+1/2}$, we have
	\begin{align}
		f^{\prime}_{*, j} = \frac{(\hat{f}_{*, j+1/2} - \hat{f}_{*, j-1/2})}{\delta}.
		\label{eq:A.2}
	\end{align}
	
	According to the definition of $\Phi(\varphi_n)$ and $\varphi_n$, we have
	\begin{align}
		\left|\Phi(\varphi_n) - \varphi_n\right| & = 
		\left|-i\frac{ \mathcal{F}[f^\prime_j](\varphi_n) \delta}{ a \mathcal{F}[u](w_n)} + i \frac{ \mathcal{F}[f^\prime](\varphi_n) \delta}{ a \mathcal{F}[u](w_n)}\right| \notag \\
		& = \frac{\delta}{a} \left| \frac{1}{N} \sum_{j = 0}^{N-1} f^\prime_{j} e^{-i\varphi_n j} - \frac{1}{N} \sum_{j = 0}^{N-1} f^{\prime}(x_j) e^{-i\varphi_n j}\right| \notag \\
		& = \frac{\delta}{a} \left| \frac{1}{N} \sum_{j = 0}^{N-1} \left[(f^\prime_{R, j} - f^\prime_R(x_j)) + i(f^\prime_{I, j} - f^\prime_I(x_j))\right](\cos(j\varphi_n) - i\sin(j\varphi_n)) \right| \notag \\
		& \le \frac{\delta}{a} \left( \frac{4}{N} \sum_{j = 0}^{N-1} \left((f^\prime_{R, j} - f^\prime_R(x_j))^2 + (f^\prime_{I, j} - f^\prime_I(x_j))^2\right) \right) ^{1/2}. \label{eq:A.3}
	\end{align}
	From the Eq. (\ref{eq:A.1}) and (\ref{eq:A.2}), we have
	\begin{subequations}
		\begin{align}
			(f^\prime_{R, j} - f^\prime_R(x_j))^2 
			\le  & \frac{2}{\delta^2} \left[ (\hat{f}_{R, j+1/2} - h_{R, j+1/2})^2 + (\hat{f}_{R, j-1/2} - h_{R, j-1/2})^2 \right], \label{eq:A.6a} \\
			(f^\prime_{I, j} - f^\prime_I(x_j))^2 
			\le  & \frac{2}{\delta^2} \left[ (\hat{f}_{I, j+1/2} - h_{I, j+1/2})^2 + (\hat{f}_{I, j-1/2} - h_{I, j-1/2})^2 \right]. \label{eq:A.6b}
		\end{align}
		\label{eq:A.6}
	\end{subequations}
	Taking the Eq. (\ref{eq:A.6}) into Eq. (\ref{eq:A.3}), the spectral error is obtained as
	\begin{align}
		\left|\Phi(\varphi_n) - \varphi_n\right| \le & \frac{1}{a} \left[ \frac{8}{N} \sum_{j = 0}^{N-1} \left( (\hat{f}_{R, j+1/2} - h_{R, j+1/2})^2 + (\hat{f}_{R, j-1/2} - h_{R, j-1/2})^2 \right. \right. \notag \\
		& \left. \left. (\hat{f}_{I, j+1/2} - h_{I, j+1/2})^2 + (\hat{f}_{I, j-1/2} - h_{I, j-1/2})^2 \right) \right]^{1/2}.
		\label{eq:A.7}
	\end{align}
	
	Eq. (\ref{eq:A.7}) provides an upper bound for the spectral error of nonlinear schemes at a specified wavenumber $\varphi_n$, demonstrating that reducing the discrepancy $\hat{f}_{*, j\pm1/2} - h_{*, j\pm1/2}$ is beneficial for improving both dispersion and dissipation. This insight offers a new perspective for enhancing the spectral properties of nonlinear numerical schemes. In the following, we utilize Eq. (\ref{eq:A.7}) to improve the spectral accuracy of the optimized WENO schemes.
	
	\subsection{WENO5-NN method}
	\label{section3.4}
	
	The optimized WENO method is defined as
	\begin{align}
		\bm{w} = \bm{w}^{*} + \bm{w}^{NN}.
		\label{eq:3.2.1}
	\end{align}	
	$\bm{w}^{*} = [w_0^{*}, w_1^{*}, w_2^{*}]$ and $*$ is JS or Z. $\bm{w}^{NN} = [w_0^{NN}, w_1^{NN}, w_2^{NN}]$ is the compensation term learned by the neural network and $\bm{w}$ is the weight function of WENO5-NN. $\bm{w}^{*}$ provides a good set of weights for the stencils in both smooth and discontinuous regions, which can reduce the difficulty of directly learning $\bm{w}$.
	
	\subsubsection{WENO5-NN structure}
\begin{figure}
		\centering
		\includegraphics[scale=0.3]{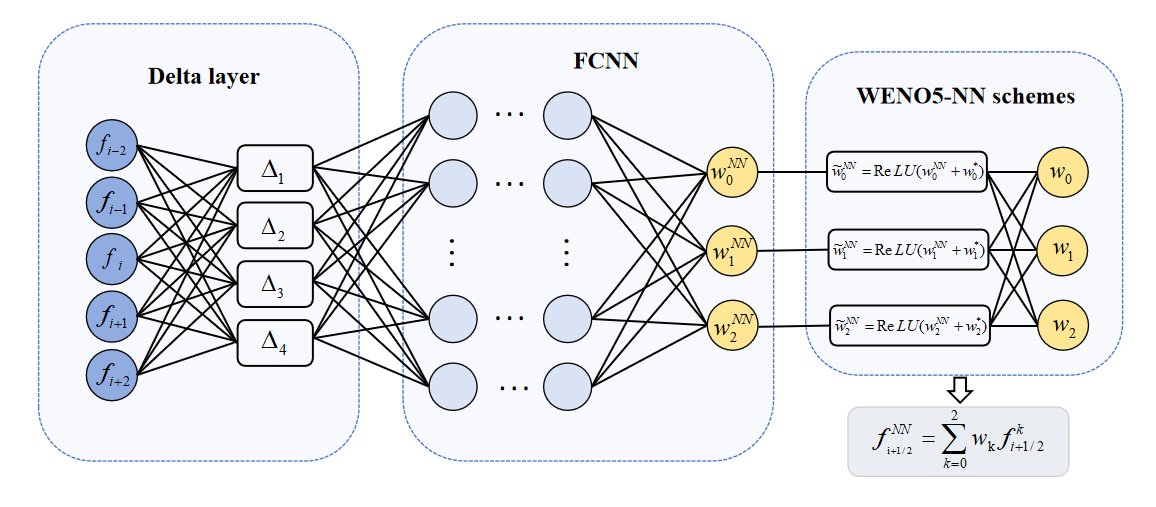}
		\caption{Structure of the neural network.}
		\label{fig:3.1}
	\end{figure}

	From section \ref{section2}, we have known that the nonlinear weight vector $\bm{w}$ is a vector function of the stencil $\{f_{i-2}, f_{i-1}, f_{i}, f_{i+1}, f_{i+2}\}$. Therefore, we denote a neural network $\mathcal{N}: \mathbb{R}^5 \rightarrow \mathbb{R}^3$, i.e., 
	\begin{align}
		\bm{w}^{NN} = \mathcal{N} (f_{i - 2}, f_{i-1}, f_{i},   f_{i+1}, f_{i+2}; \theta),
		\label{eq:3.2.2}
	\end{align}
	whose structure is displayed in Fig. \ref{fig:3.1}. 
	The input values are the fluxes $\{f_{i-2}, f_{i-1}, f_{i}$, $f_{i+1}, f_{i+2}\}$. First, the vector is passed to the Delta layer, which calculates the features from the input values
	\begin{equation}
		\begin{aligned}
			& \tilde{\Delta}_j = |f_{i-2+j} - f_{i-3+j}|, \quad \Delta_j = \tilde{\Delta}_j /   \max(\tilde{\Delta}_j,  \epsilon_1), \quad \quad j = 1, 2, 3, 4.
			\label{eq:3.2.3}
		\end{aligned}
	\end{equation}
	$\epsilon_1$ is a small value used to avoid a zero denominator. $\Delta_j$ is a good measurement of the regularity of the local flux field. Then the features $\{\Delta_1, \Delta_2, \Delta_3, \Delta_4 \}$ are passed through a multilayer perceptron. The output vector is $\bm{w}^{NN}$.
	
	To satisfy the condition that $\sum_{l=0}^{2} w_l = 1$, $w_l \ge 0$, the weights in Eq. (\ref{eq:3.2.1}) are normalized by
	\begin{align}
		w_k =  \frac{\tilde{w}_k}{\sum_{l=0}^{2}\tilde{w}_l}, \quad \tilde{w}_k = ReLU(w_k^* + w_k^{NN}), \quad k = 0, 1, 2.
		\label{eq:3.2.4}
	\end{align}	
	
	\subsubsection{Loss function}
	
	To enhance the spectral accuracy while preserving the essential non-oscillatory behavior near discontinuities, we construct four loss terms, three of which are tailored to improve specific properties.

	The first loss term $\mathcal{L}_r$ is the interface reconstruction error, defined as
	\begin{equation}
		\mathcal{L}_r = \frac{1}{N_b N_x} \sum_{s=1}^{N_b} \sum_{i=0}^{N_x} \left(\hat{f}_{i+\frac{1}{2}}^{NN, s} -  h_{i+\frac{1}{2}}^{s}\right)^2. \label{eq:3.2.6}	\\
	\end{equation}
	$N_x$ is the mesh grid number, and $N_b$ is the mini-batch size of the training dataset. $\mathcal{L}_r$ quantifies the discrepancy between the reconstructed numerical flux $\hat{f}_{i+\frac{1}{2}}^{NN, s}$ predicted by the neural network and the reference flux $h_{i+\frac{1}{2}}^s$. By minimizing this loss over a large set of stencils generated from sine waves with different wavenumbers, the spectral accuracy of WENO5-NN can be improved.
	
	The second loss $\mathcal{L}_{TVD}$ is the TVD penalization. We apply the single-step TVD condition and transform the inequality constraint Eq. (\ref{eq:3.2.7a}) into a soft-constraint function as 
	\begin{align}
		\mathcal{L}_{TVD} = \frac{1}{N_b}\sum_{s = 1}^{N_b} \max \left(  \text{TV}(u_{\Delta t}^{NN, s}) - \text{TV}(u^{s}), 0 \right)^2,
		\label{eq:3.2.7}
	\end{align}
	where $u_{\Delta t}^{NN, s}$ is computed based on the initial condition $u^{s}$, the predicted interface flux $\hat{f}_{i+\frac{1}{2}}^{NN, s}$ in space, and the forward Euler method for time integration. The TVD penalty effectively inhibits the non-physical oscillations caused by the WENO5-NN scheme near discontinuities.
	
	The third loss $\mathcal{L}_{diss}$ is the anti-dissipation penalization. Though the new numerical scheme demonstrates favorable spectral properties under the loss function Eq. (\ref{eq:3.2.6}), it inevitably produces anti-dissipation, which can compromise the numerical stability. Therefore, it is necessary to suppress anti-dissipation.
	
	To restrain the generation of the anti-dissipation in the training process, we introduce the imaginary part of the modified wavenumber $\Phi(\varphi_n)$ to the loss function as
	\begin{align}
		\mathcal{L}_{diss} = \frac{1}{N/2}\sum_{n = 0}^{N/2} \max \left(  \text{Im} (\Phi(\varphi_n)), 0 \right)^2,
		\label{eq:3.2.9}
	\end{align}
	where the spatial derivative for $\Phi(\varphi_n)$ is calculated by the training WENO5-NN scheme. The loss function $\mathcal{L}_{diss}$ guides the optimization of the WENO5-NN scheme toward positive or non-dissipative behavior by penalizing anti-dissipative effects, that is, $\text{Im} (\Phi(\varphi)) \ge 0$.
	
	The last loss $||\bm{W}||_2^2$ is an $l_2$-regularization loss to prevent overfitting, where $\bm{W}: = \{\bm{W}^{[l]}\}_{l=1}^{L}$ in Eq. (\ref{eq:3.1.1}) is the set of weight matrices of the neural network.
	
	The total loss is defined as a combination of the interface reconstruction error, TVD penalization, anti-dissipation penalization and $l_2$-regularization:
	\begin{equation}
		\mathcal{L} = \mathcal{L}_r + \lambda_{TVD}  \mathcal{L}_{TVD} + \lambda_{diss} \mathcal{L}_{diss} + \lambda_W ||\bm{W}||_2^2. \label{eq:3.2.5}
	\end{equation}	
	$\lambda_{TVD}$, $\lambda_{diss}$ and $\lambda_{W}$ are positive hyperparameters. 
	
	\section{Numerical examples}
	\label{section4}
	\FloatBarrier
	In this section, several benchmark examples of the Euler equations are used to illustrate the potential of the proposed methods for high-resolution simulations. Throughout the following tests, we use the global Lax--Friedrichs method for flux splitting and the Roe approximation for the characteristic decomposition at the cell faces. The third-order TVD Runge--Kutta method \cite{ref9} is adopted to integrate the equations in time with $\text{CFL} = 0.4$. The numerical results are compared with WENO5-JS and WENO5-Z.
	
	The neural network consists of an input layer, three hidden layers with $30$ nodes, and an output layer with $3$ nodes. The rectified linear unit (ReLU) activation function is used for hidden layers. During training, the learning rate is $10^{-3}$ and decays exponentially. The mini-batch $N_b$ is $800$ and the hyperparameters $\lambda_{W}$ is $10^{-8}$. $\epsilon_1 = 10^{-30}$ is used in the Eq. (\ref{eq:3.2.3}). 
	
	\subsection{Training dataset}
	\label{section4.1}
	
	\begin{table}[h]
		\centering
		\caption{Training dataset.}
		\begin{tabular}{ccc}
			\toprule
			Function $f(x)$		& Parameters	& Number of samples	\\
			\hline
			$\tanh(ax)$		&	$|a| \sim \mathcal{U}(50, 100)$	&	2000	\\
			$\sin(b \pi x + \phi)$	& $b \sim \mathcal{U}(1, 18)$, $\phi  \sim \mathcal{U}(0, 2\pi)$ 	&	1000	\\
			$\sum_{k=0}^{5} c_k x^k$ & $c_k \sim \mathcal{U}(-\frac{1}{k+1}, \frac{1}{k+1})$	&	1000	\\
			\bottomrule
		\end{tabular}
		\label{tab:3.1}
	\end{table}

	The training dataset is composed of canonical functions that resemble the local features of solutions to hyperbolic conservation laws, as detailed in Table \ref{tab:3.1}, where $\mathcal{U}$ represents a uniform distribution.
	We use three classes of functions: $\tanh$ functions, sinusoidal functions, and polynomials up to degree $5$. All the functions are calculated on the domain $[-1, 1]$, discretized with $N_x = 100$. 
	The $\tanh$ functions with $|a|>50$ exhibit sharp variations at $x=0$ and can be used to mimic discontinuity. 
	The sinusoidal function can be used to improve the spectral properties of the WENO5-JS/Z-NN schemes. The wavenumber range has an impact on the numerical spectrum. 
	Here, we chose $b \in [1, 18]$, which corresponds to scaled wavenumber $\varphi \in [0.02, 1.1310]$. 
	We apply a high-order numerical method to approximate the exact interface flux $h_{i+\frac{1}{2}}$ by expanding the flux function $f(u)$ at $x_{i+1/2}$. Following \cite{ref6}, the formulation is 
	\begin{align}
		\hat{f}_{i+\frac{1}{2}} = f_{i+\frac{1}{2}} -  \frac{1}{24}\Delta x^2 	f_{xx}\bigg|_{i+\frac{1}{2}} + \frac{7}{5760}\Delta x^4 f_{xxxx}\bigg|_{i+\frac{1}{2}} + O(\Delta x^6).
		\label{eq:3.4.1}
	\end{align}

	\subsection{Validity test of the penalization}		
	\label{section4.2}
	To verify the validity of the regularization terms $\mathcal{L}_{TVD}$ and $\mathcal{L}_{diss}$, we compare WENO5-Z-NN methods on discontinuous and high-frequency problems under different hyperparameter settings.
	
	\noindent \textbf{TVD penalization} 
	
\begin{figure}[htbp]
		\centering
		\subfloat[WENO5-Z-NN\label{fig:4.1.1a}]
		{\includegraphics[width=0.45\textwidth]{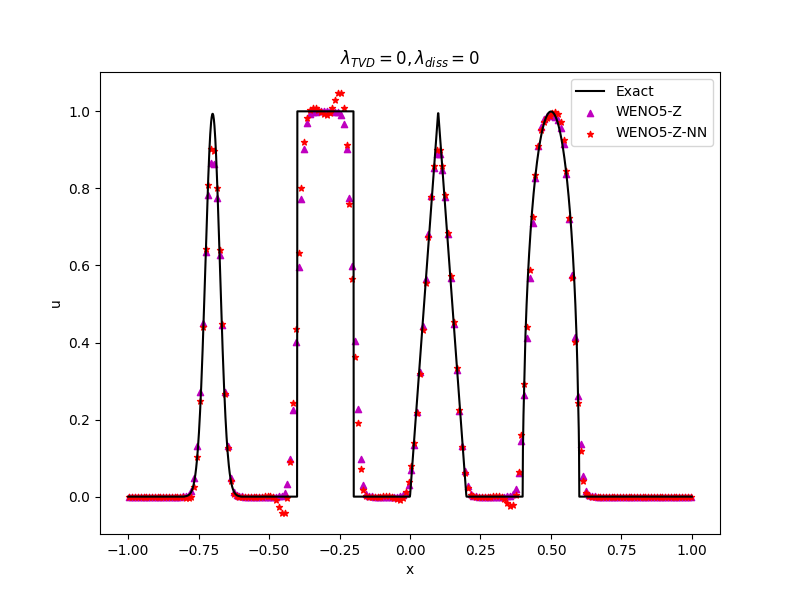}}
		\hspace{0.1cm}
		\subfloat[WENO5-Z-NN]
		{\includegraphics[width=0.45\textwidth]{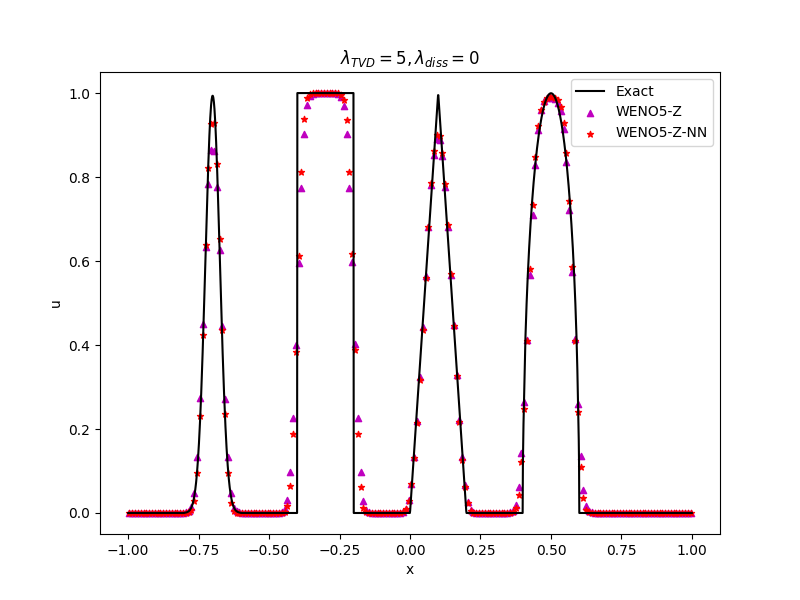}\label{fig:4.1.1b}}
		\caption{Linear advection equation calculated by WENO5-Z-NN at $t = 4.0$ with $N_x = 200$.}
		\label{fig:4.1.1}
	\end{figure}

	WENO5-Z-NN with hyperparameters $(\lambda_{TVD}, \lambda_{diss}) = (0, 0)$ and $(5, 0)$ are compared to demonstrate the effect of the TVD penalization. We consider the 1D linear advection problem with a discontinuous initial condition consisting of a Gaussian, a square, a triangle, and a semi-ellipse, given by
	
	\begin{align}
		u(x, t = 0) = \begin{cases}
			\frac{1}{6}\left( G(x, \beta, z - \delta) +  G(x, \beta, z + \delta) + 4 G(x, \beta, z) \right) \quad & -0.8 \le x < -0.6, \\
			1  \quad &  -0.4 \le x < -0.2 ,	\\
			1 - |10 (x - 0.1)|  \quad &  0.0 \le x <  0.2, \\
			\frac{1}{6}\left( F(x, \alpha, a - \delta)  + F(x, \alpha, a + \delta) + 4 F(x, \alpha, a) \right)  \quad &  0.4 \le x < 0.6, \\
			0   \quad &  \text{otherwise},
		\end{cases}
		\notag
	\end{align}
	where $G(x, \beta, z) = e^{-\beta (x - z)^2}$, $F(x, \alpha, a) = \sqrt{\max(1 - \alpha^2(x - a)^2, 0)}$. The constants are $z = -0.7$, $\delta = 0.005$, $\beta = \frac{\log 2}{36 \delta^2}$, $a = 0.5$, and $\alpha = 10$. The domain  is $[-1, 1]$ and the resolution is $N_x = 200$. The periodic boundary conditions are applied. The results at $t = 4$ are shown in Fig. \ref{fig:4.1.1}. 
	Both schemes successfully resolve the solution profile. Nonetheless, WENO5-Z-NN with $\lambda_{TVD} = 0$ exhibits evident oscillations around the square and semi-elliptic waves. The method incorporating TVD penalization with $\lambda_{TVD} = 5$ effectively mitigates such oscillations.
	
	\noindent \textbf{Anti-dissipation penalization} 
	
\begin{figure}[htbp]
		\centering
		\subfloat[WENO5-Z-NN]
		{\includegraphics[width=0.45\textwidth]{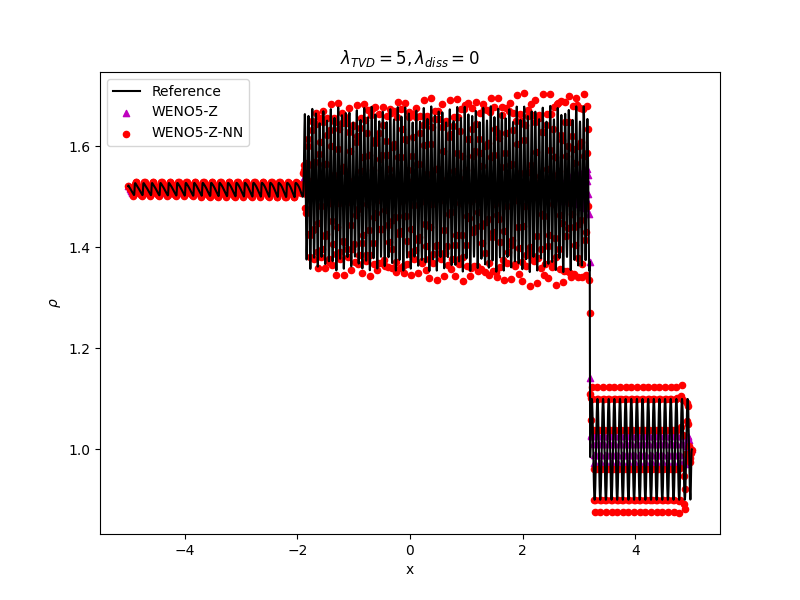}\label{fig:4.1.2a}}
		\hspace{0.1cm}
		\subfloat[WENO5-Z-NN]
		{\includegraphics[width=0.45\textwidth]{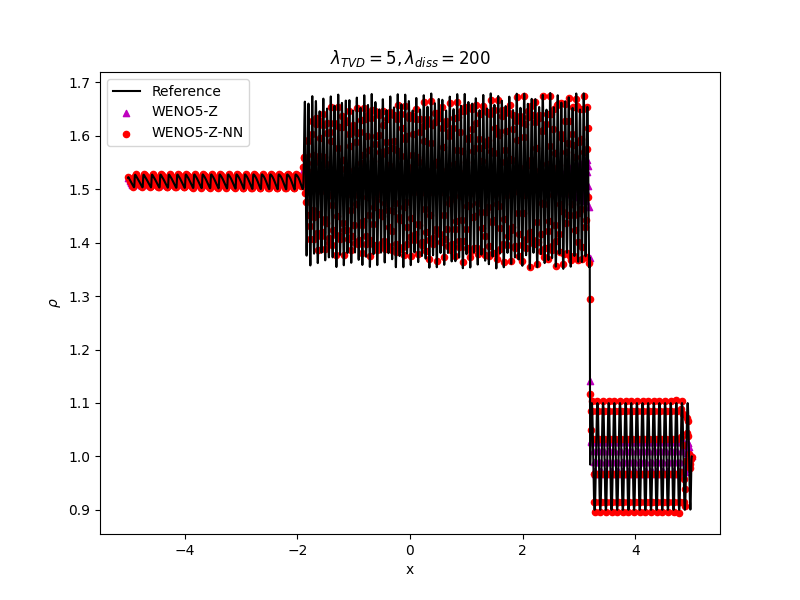}\label{fig:4.1.2b}}
		\caption{Density distribution of the Titarev and Toro problem calculated by WENO5-Z-NN at $t = 5.0$ with $N_x = 1000$.}
		\label{fig:4.1.2}
	\end{figure}

	To explore the effect of the anti-dissipation penalization, we test the WENO5-Z-NN with hyperparameters $(\lambda_{TVD}, \lambda_{diss}) = (5, 0)$ and $(5, 200)$ on the Titarev and Toro problem \cite{ref45}. The initial condition on the computational domain $[-5, 5]$ is 
	\begin{equation}
		\begin{aligned}
			(\rho, u, p) =
			\begin{cases}
				(1.515695, 0.523346, 1.805), & -5 \le x < -4.5, \\
				(1 + 0.1 \sin(20 \pi x), 0, 1.0), & -4.5 \le x \le 5.	
			\end{cases}
			\label{eq:4.1.1}
		\end{aligned}
	\end{equation}
	Non-reflection boundary conditions are employed on the left and right boundaries. The final simulation time is $t = 5$. Fig. \ref{fig:4.1.2} compares the density distribution obtained by WENO5-Z-NN on a grid $N_x = 1000$, where the reference solution is calculated by the WENO5-JS scheme with $N_x = 5000$. The results show that these two schemes can capture high-frequency waves. However, WENO5-Z-NN with $\lambda_{diss} = 0$ amplifies the amplitude of the waves within the domain $[-2, 5]$, resulting in noticeable overshoots relative to the reference solution. In contrast, the WENO5-Z-NN scheme with $\lambda_{diss} = 200$ does not exhibit such issues, demonstrating the effectiveness of the anti-dissipation penalization.
	
	\subsection{Spectral properties and the distribution of weights}
	\label{section4.3}
\begin{figure}[htbp]
		\centering
			\subfloat[Dispersion of the JS schemes]
			{\includegraphics[width=0.4\textwidth]{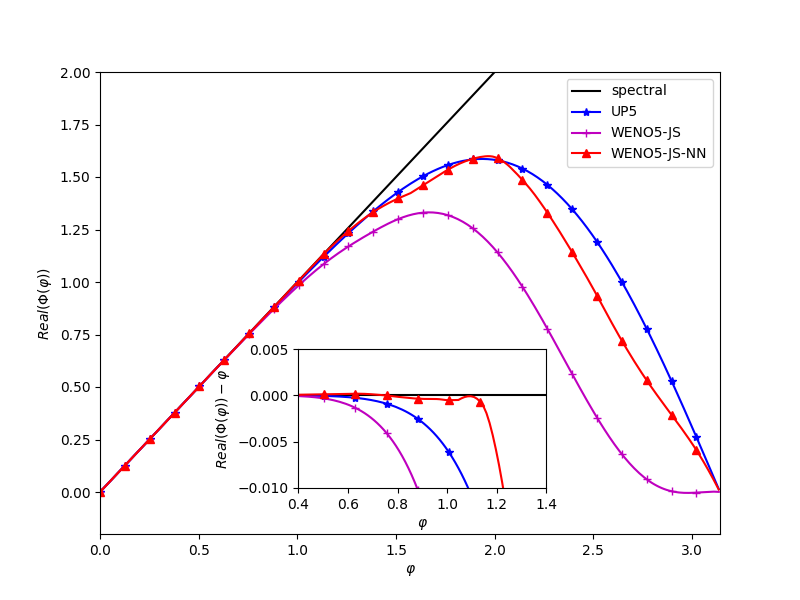}}
			\subfloat[Dissipation of the JS schemes]
			{\includegraphics[width=0.4\textwidth]{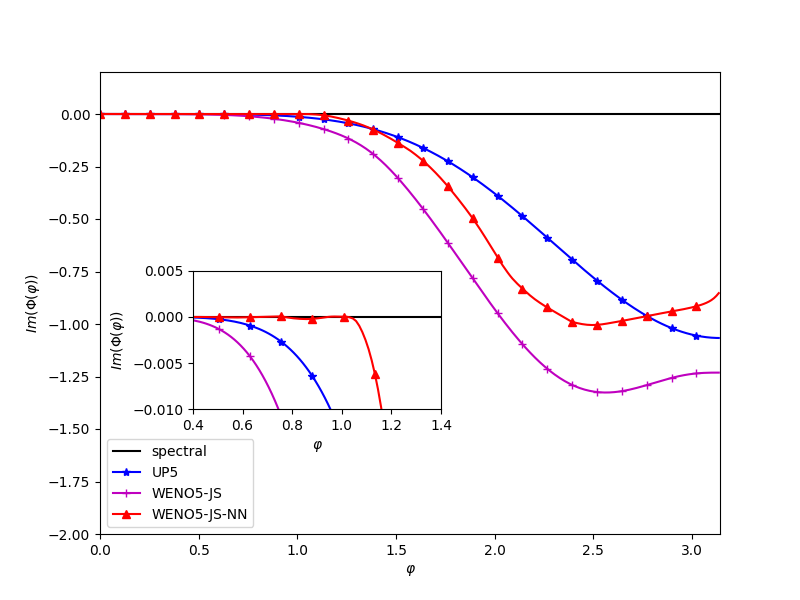}}
			\\
			\subfloat[Dispersion of the Z schemes]
			{\includegraphics[width=0.4\textwidth]{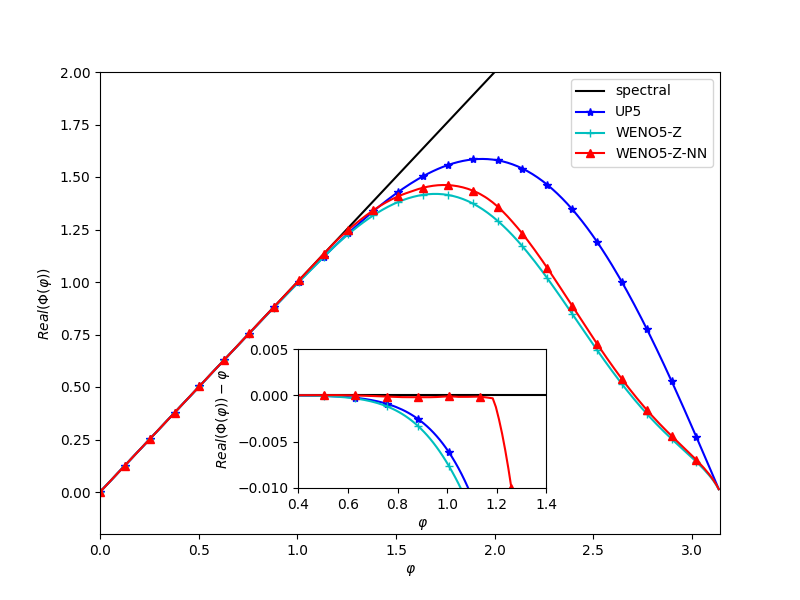}}
			\subfloat[Dissipation of the Z schemes]
			{\includegraphics[width=0.4\textwidth]{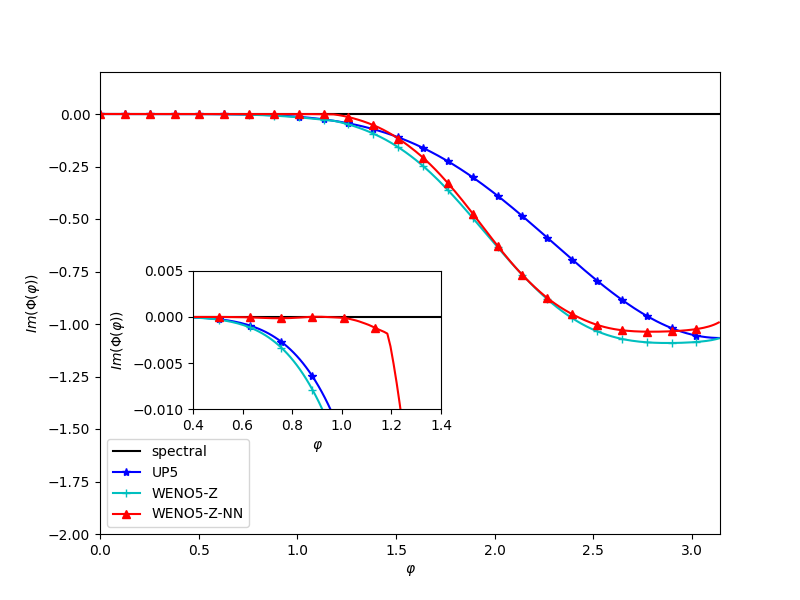}}
		\caption{Approximate dispersion and dissipation properties of the UP5, WENO5-JS/Z and WENO5-JS/Z-NN schemes.}
		\label{fig:4.2.1}
	\end{figure}

	In the subsequent simulations, the hyperparameters ($\lambda_{TVD}, \lambda_{diss}$) of WENO5-JS-NN are set to $(80, 700)$, while those of WENO5-Z-NN are set to $(5, 200)$.
	Fig. \ref{fig:4.2.1} shows the dispersion and dissipation of the UP5, WENO5-JS/Z, WENO5-JS/Z-NN schemes. All the schemes have good spectral properties for small wavenumbers. The WENO5-JS scheme is more dissipative than WENO5-Z and has more room for optimization.
	For the medium wavenumbers ($ 0.4 \le \varphi \le 1.4$), the numerical spectrum of the WENO5-JS/Z-NN schemes is close to the exact spectrum, which indicates that the new schemes have lower dispersion and dissipation.
	Furthermore, the enlarged views in Fig.~\ref{fig:4.2.1} demonstrate that our schemes exhibit enhanced spectral accuracy within the learned wavenumber interval ($\varphi \in [0.02, 1.1310]$), outperforming the UP5 scheme.
	
\begin{figure}[htbp]
		\begin{tabular}{c@{\hspace{-3pt}}c@{\hspace{-3pt}}c@{\hspace{-3pt}}c}
			\subfloat[JS schemes]{\includegraphics[width=0.25\textwidth]{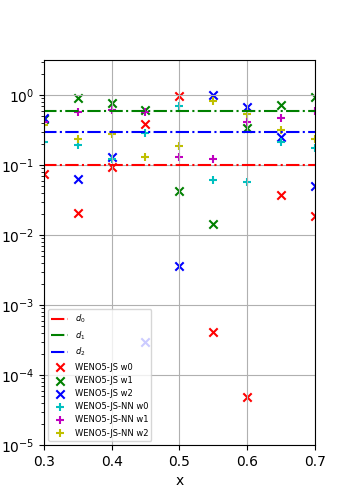}} & \subfloat[Z schemes]{\includegraphics[width=0.25\textwidth]{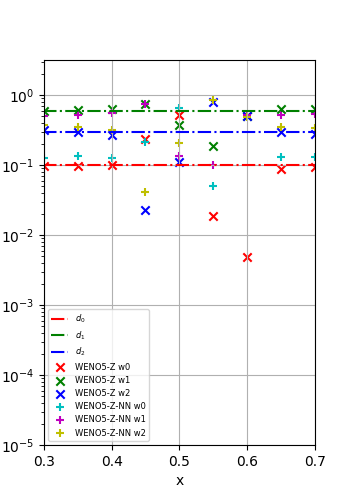}} &
			\subfloat[JS schemes]{\includegraphics[width=0.25\textwidth]{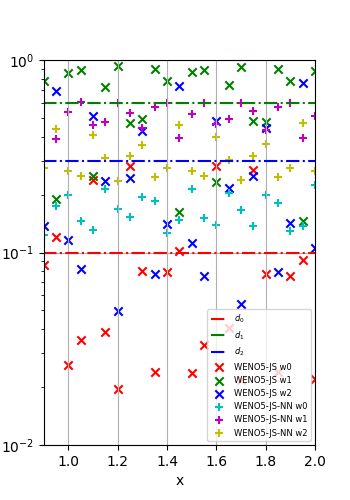}}	&
			\subfloat[Z schemes]{\includegraphics[width=0.25\textwidth]{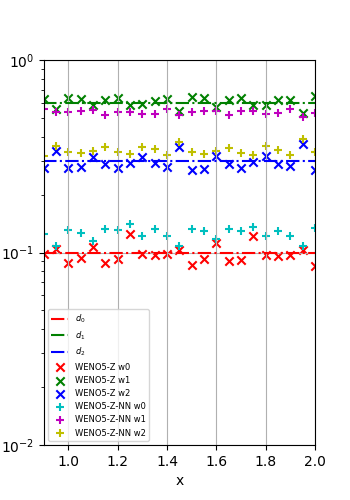}}
		\end{tabular}
		\caption{Distribution of the ideal weights $d_k$ and nonlinear weights $w_k$ for different schemes. (a) and (b) are near the discontinuity. (c) and (d) are in the smooth region.}
		\label{fig:4.2.3}
	\end{figure}

	To analyze the distribution of the nonlinear weights for WENO5-JS/Z-NN, we consider the function as 
	\begin{align}
		f(x) = \begin{cases}
			\frac{2}{3} \sin(6 \pi x) + \frac{1}{4} \sin(1.6 \pi x), \quad & 0 \le x < 0.5,	\\	\notag	
			\frac{2}{3} \sin(6 \pi x) + \frac{1}{4} \sin(1.6 \pi x) + 3,  \quad & 0.5\le x \le 2,
		\end{cases}
	\end{align}
	which contains a discontinuity at $x = 0.5$. The shape of the smooth wave is similar to the high-frequency component observed within the density distribution of the Shu–Osher problem. Fig. \ref{fig:4.2.3} exhibits the distribution of the nonlinear weights for these schemes near the discontinuity and in the smooth region, respectively. 
	WENO5-JS/Z-NN assign small weights to stencils containing the discontinuity. In particular, we find that the WENO5-Z-NN scheme attempts to increase the contribution of the stencil that is less smooth.
	For the smooth wave, the weights of WENO5-JS/Z fluctuate around the ideal linear weights, and the amplitude of WENO5-JS is larger than that of WENO5-Z. While the WENO5-JS/Z-NN schemes seem to modify background linear weights by increasing the weights $w_0$ and $w_2$ while decreasing the weight $w_1$, aligning with the results adopted in linear spectral optimization \cite{ref53}.

	\subsection{Numerical examples for the 1D Euler equations}
	\label{section4.4}
	Consider the 1D Euler equations 
	\begin{align}
		\frac{\partial \bm{U}}{\partial t}
		+
		\frac{\partial \bm{F}}{\partial x}
		= 0,
		\label{Eq:4.3}
	\end{align}
	where $\bm{U} = [\rho, \rho u, E]$ is a conservation vector and $\bm{F} = [\rho u, \rho u^2 + p, u(E+p)]$ is a flux vector. $\rho$, $u$, $E$, and $p$ are density, velocity, total energy, and pressure, respectively. The system represents the conservation of mass, momentum, and energy. The total energy for an ideal polytropic gas is defined as 
	\begin{align*}
		E = \frac{p}{\gamma - 1} + \frac{1}{2} \rho u^2,
	\end{align*}
	where $\gamma$ is the specific heat ratio and is set to $1.4$.
	
	\noindent \textbf{Example 1.} Lax problem.
	
\begin{figure}[htbp]
		\centering
		\subfloat[WENO5-JS and WENO5-JS-NN]
		{\includegraphics[width=0.45\textwidth]{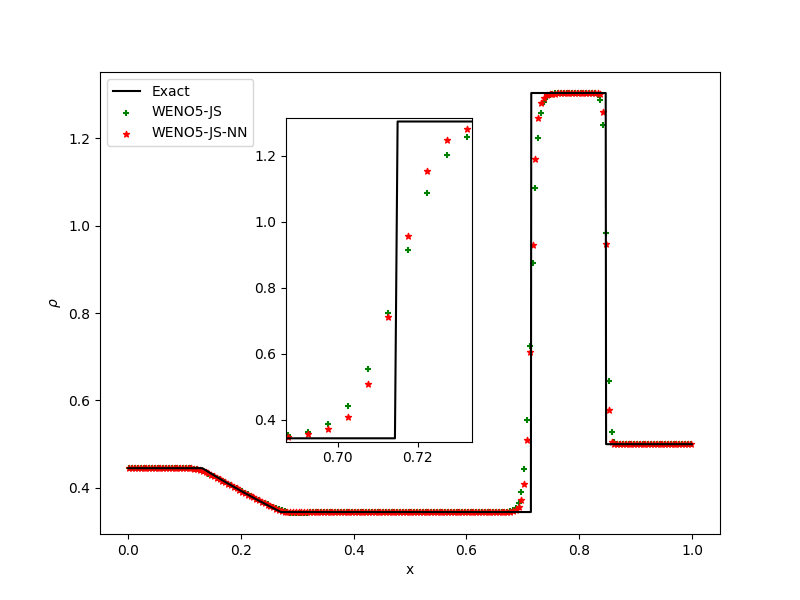}\label{fig:4.3.2a}}
		\hspace{0.1cm}
		\subfloat[WENO5-Z and WENO5-Z-NN]
		{\includegraphics[width=0.45\textwidth]{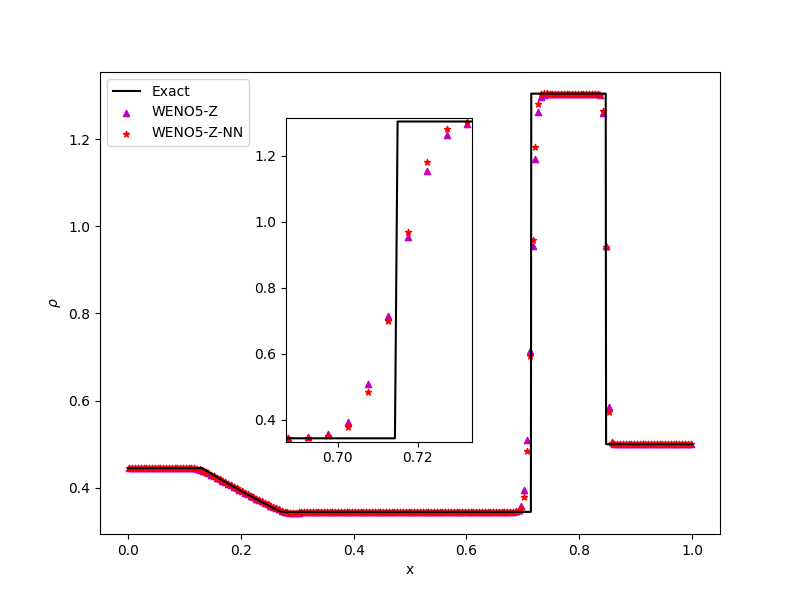}\label{fig:4.3.2b}}
		\caption{Density distribution of the Lax problem at $t = 0.14$ with $N_x = 200$.}
		\label{fig:4.3.1}
	\end{figure}

	The initial condition of the Lax problem \cite{ref47} is 
	\begin{equation}
		\begin{aligned}
			(\rho, u, p) =
			\begin{cases}
				(0.445, 0.698, 3.528), & 0 \le x < 0.5, 	\\
				(0.5, 0, 0.5710), & 0.5 \le x <1.0.	
			\end{cases}
			\notag	
		\end{aligned}
	\end{equation}
	Non-reflection boundary conditions are employed. The numerical simulations are advanced till $t = 0.14$, and the grid is $N_x = 200$. The density solutions are compared in Fig. \ref{fig:4.3.1}. All the schemes capture the contact discontinuity and the right-moving shock, where the results of WENO5-JS/Z-NN are much sharper than those of WENO5-JS/Z.		
	
	\noindent \textbf{Example 2.} Woodward--Colella blast waves.
	
\begin{figure}[htbp]
		\centering
		\subfloat[WENO5-JS and WENO5-JS-NN]
		{\includegraphics[width=0.45\textwidth]{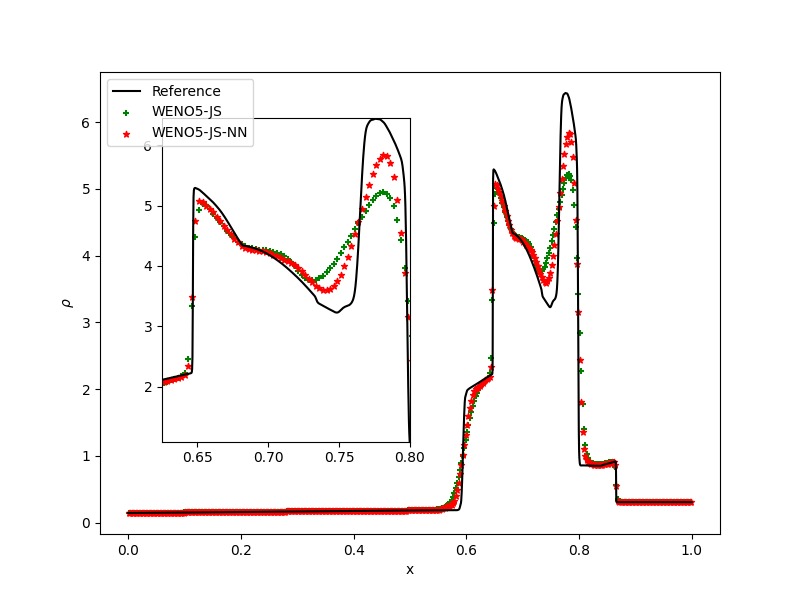}\label{fig:4.3.1a}}
		\hspace{0.1cm}
		\subfloat[WENO5-Z and WENO5-Z-NN]
		{\includegraphics[width=0.45\textwidth]{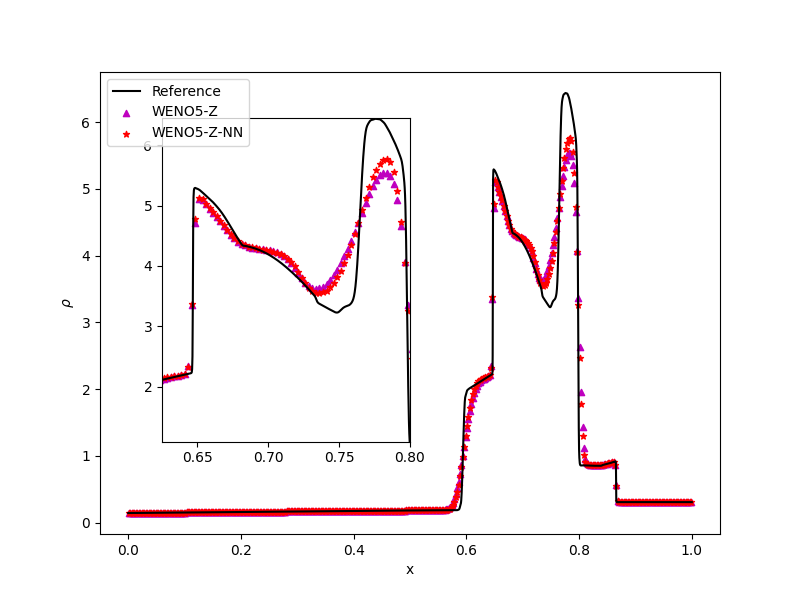}\label{fig:4.3.1b}}
		\caption{Density distribution of the Woodward--Colella blast wave problem at $t = 0.038$ with $N_x = 400$.}
		\label{fig:4.3.2}
	\end{figure}

	The Woodward--Colella blast wave problem \cite{ref46} is a high Mach number test used to evaluate the ability of numerical methods to resolve strong shocks and their interactions. The initial condition is 
	\begin{equation}
		\begin{aligned}
			(\rho, u, p) =
			\begin{cases}
				(1, 0, 1000), & 0 \le x < 0.1, \\
				(1, 0, 0.01), & 0.1 \le x <0.9,	\\
				(1, 0, 100), & 0.9 \le x <1.0.	\\
			\end{cases}
			\notag
		\end{aligned}
	\end{equation}
	Reflection boundary conditions are applied on the left and right. Fig. \ref{fig:4.3.2} compares the density distribution computed by WENO5-JS/Z and WENO5-JS/Z-NN with the resolution $N_x = 400$, where the reference solution is calculated by the WENO5-JS scheme with $N_x = 5000$. The final simulation time is $t = 0.038$. Our schemes exhibit less dissipation near the valley at $x = 0.75$ and the right peak at $x = 0.78$, while providing more accurate results of the contact discontinuity and the post-impact structure.
	
	\noindent \textbf{Example 3.}  Shu--Osher problem.
	
\begin{figure}[htb]
		\centering
		\subfloat[WENO5-JS and WENO5-JS-NN]
		{\includegraphics[width=0.45\textwidth]{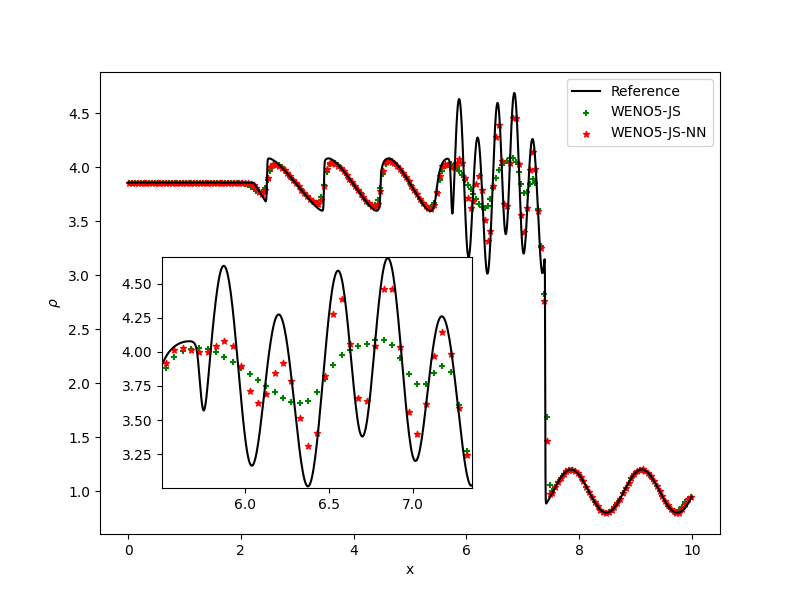}\label{fig:4.3.3a}}
		\hspace{0.1cm}
		\subfloat[WENO5-Z and WENO5-Z-NN]
		{\includegraphics[width=0.45\textwidth]{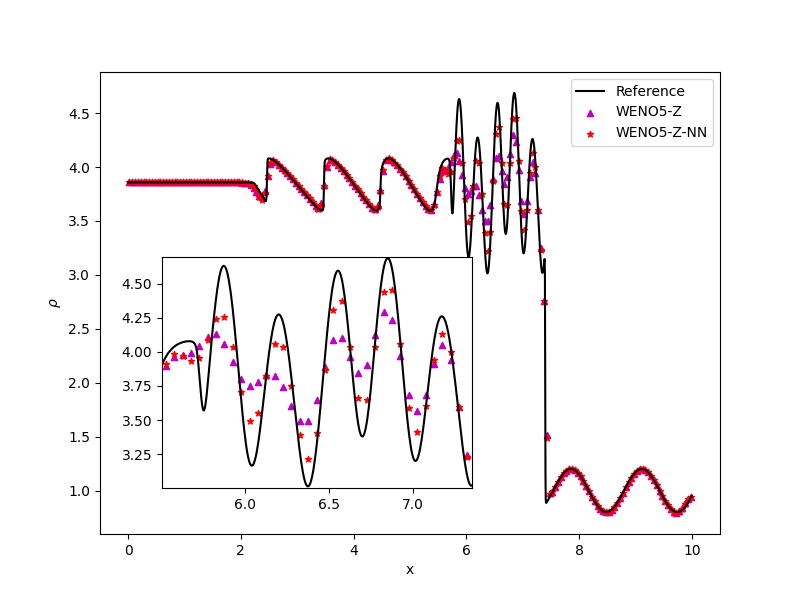}\label{fig:4.3.3b}}
		\caption{Density distribution of the Shu--Osher problem at $t = 1.8$ with $N_x = 200$.}
		\label{fig:4.3.3}
	\end{figure}

	The Shu--Osher problem \cite{ref6} contains a Mach 3 shock interacting with a perturbed density field and is widely used to evaluate the stability and accuracy of the investigated scheme. The initial condition is 
	\begin{equation}
		\begin{aligned}
			(\rho, u, p) =
			\begin{cases}
				(3.857, 2.629, 10.333), & 0 \le x < 1, \\
				(1 + 0.2 \sin(5x), 0, 1.0), & 1 \le x < 10.	
			\end{cases}
			\notag	
		\end{aligned}
	\end{equation}
	Non-reflection boundary conditions are used. The final simulation time is $t = 1.8$. The resolution is $N_x = 200$. The density distribution is shown in Fig. \ref{fig:4.3.3}, where the reference solution is based on the WENO5-JS scheme with $N_x = 2000$. Compared with WENO5-JS, WENO5-JS-NN identifies more waves. The results of the WENO5-JS/Z-NN schemes are closer to the reference solution in the high-frequency region, indicating that the new schemes have lower dissipation.
	
	\noindent \textbf{Example 4.} Titarev and Toro problem.
	
	The Titarev and Toro problem is an upgraded version of the Shu--Osher problem, which describes the Mach 1.1 interaction with a high-frequency entropy sine wave. The initial condition is Eq. (\ref{eq:4.1.1}). The WENO5-JS and WENO5-Z schemes exhibit limited accuracy in resolving high-frequency wave components, while our schemes substantially improve resolution in this regime.
	
	\begin{figure}[htb]
		\centering
		\subfloat[WENO5-JS and WENO5-JS-NN]
		{\includegraphics[width=0.45\textwidth]{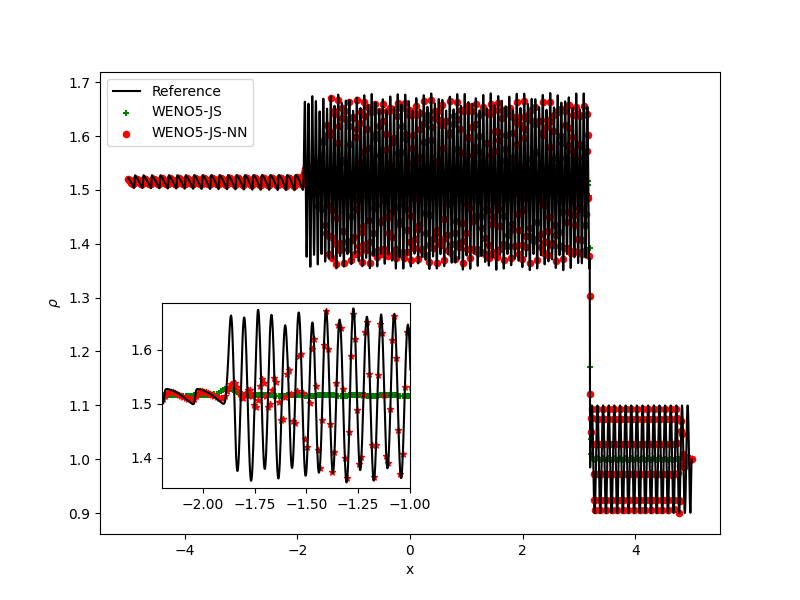}\label{fig:4.3.4a}}
		\hspace{0.1cm}
		\subfloat[WENO5-Z and WENO5-Z-NN]
		{\includegraphics[width=0.45\textwidth]{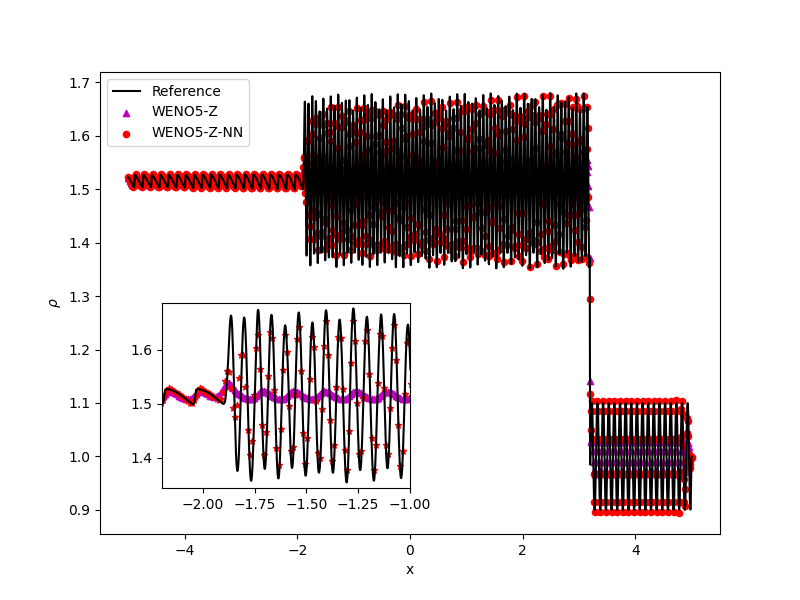}\label{fig:4.3.4b}}
		\hspace{0.1cm}
		\caption{Density distribution of the Titarev and Toro problem at $t = 5$ with $N_x = 1000$.}
		\label{fig:4.3.4}
	\end{figure}
	
	\subsection{Numerical examples for the 2D Euler equations}
	\label{section4.5}
	Consider the 2D Euler equations
	\begin{align}
		\frac{\partial \bm{U}}{\partial t}
		+
		\frac{\partial \bm{F}}{\partial x}
		+
		\frac{\partial \bm{G}}{\partial y}
		= 0,
		\label{eq:4.4.1}
	\end{align}
	where \(\bm{U} = [\rho, \rho u, \rho v, E]^T\) is the vector of conserved variables, and \(\bm{F} = [\rho u, \rho u^2 + p, \rho u v, u(E+p)]^T\) and \(\bm{G} = [\rho v, \rho u v, \rho v^2 + p, v(E+p)]^T\) are the flux vectors in the \(x\)- and \(y\)-directions, respectively.
	The total energy for an ideal polytropic gas is
	\begin{align*}
		E = \frac{p}{\gamma - 1} + \frac{1}{2} \rho  (u^2 + v^2).
	\end{align*}
	
	\noindent \textbf{Example 5.} Two-dimensional Riemann problem.
	
\begin{figure}[htb]
		\centering
		\begin{tabular}{cc}
			\subfloat[WENO5-JS]{\includegraphics[width=0.35\textwidth]{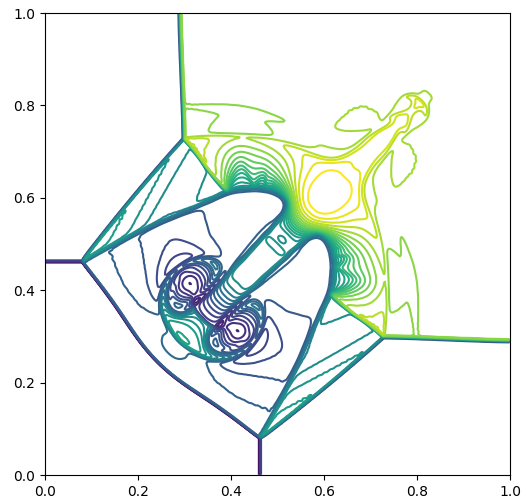}} & \subfloat[WENO5-Z]{\includegraphics[width=0.35\textwidth]{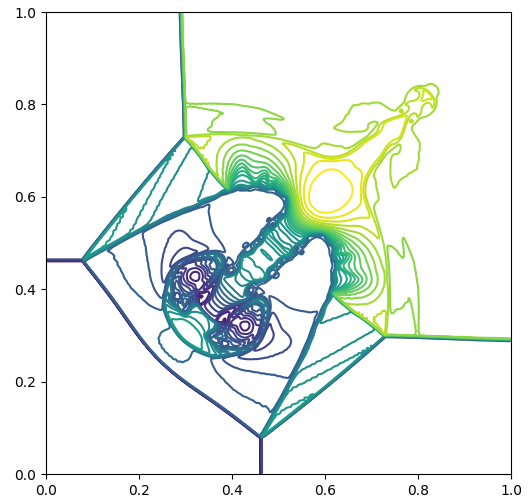}} \\
			\subfloat[WENO5-JS-NN]{\includegraphics[width=0.35\textwidth]{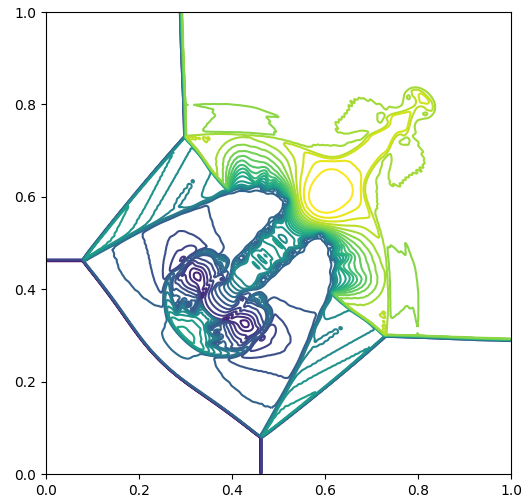}} & \subfloat[WENO5-Z-NN ]{\includegraphics[width=0.35\textwidth]{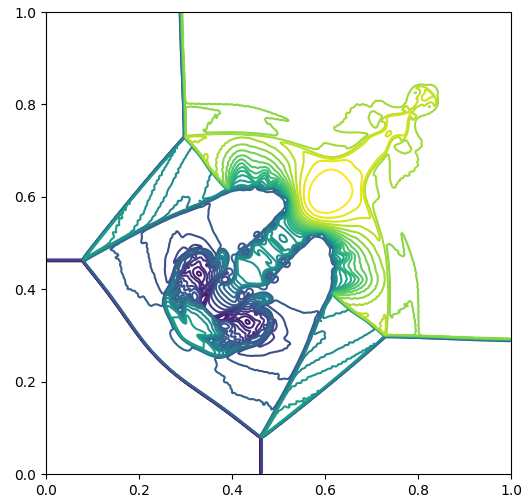}} \\
		\end{tabular}
		\caption{Density contours of the two-dimensional Riemann problem at $t = 0.8$ with $500 \times 500$ grid points. This figure is drawn with $30$ density contours between $0.2$ and $1.7$.}
		\label{fig:4.4.1}
	\end{figure}

	This test is the benchmark case 3 of the 2D Riemann problems in \cite{ref48}. The initial condition on the domain $[0, 1] \times [0, 1]$ is given as
	\begin{align}
		(\rho, u, v, p) = 
		\begin{cases}
			(1.5, 0, 0, 1.5), & 0.8 \leq x \leq 1, 0.8 \leq y \leq 1, \\
			(0.5323, 1.206, 0, 0.3), & 0 \leq x < 0.8, 0.8 \leq y \leq 1, \\
			(0.138, 1.206, 1.206, 0.029), & 0 \leq x  < 0.8, 0 \leq y < 0.8, \\
			(0.5323, 0, 1.206, 0.3), & 0.8 \leq x \leq 1, 0 \leq y < 0.8.
		\end{cases}
		\notag	
	\end{align}
	The specific heat ratio is $\gamma = 1.4$. The density contours at $t = 0.8$ calculated by WENO5-JS/Z and WENO5-JS/Z-NN with $500 \times 500$ grid points are shown in Fig. \ref{fig:4.4.1}. Each scheme can capture the reflected shock waves and contact discontinuities. WENO5-JS/Z-NN can resolve more small-scale structures near the slip line than WENO5-JS/Z. Moreover, the vortex structures computed by WENO5-Z-NN are clearer than those by WENO5-JS-NN.
	
	\noindent \textbf{Example 7.} Rayleigh--Taylor instability.
	
\begin{figure}[htbp]
		\centering
		\begin{tabular}{cccc}
			\subfloat[WENO5-JS]{\includegraphics[width=0.20\textwidth]{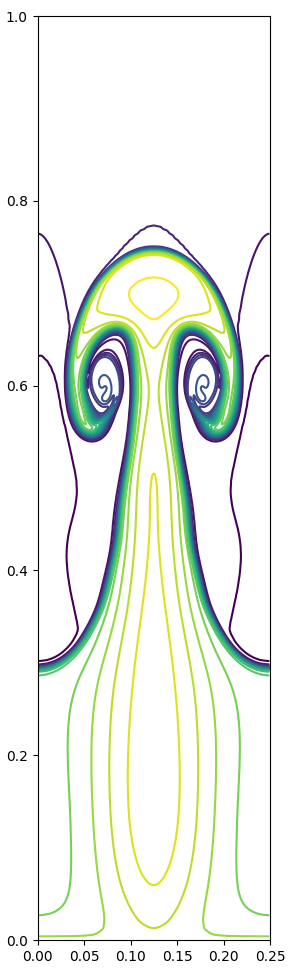}} & \subfloat[WENO5-JS-NN]{\includegraphics[width=0.20\textwidth]{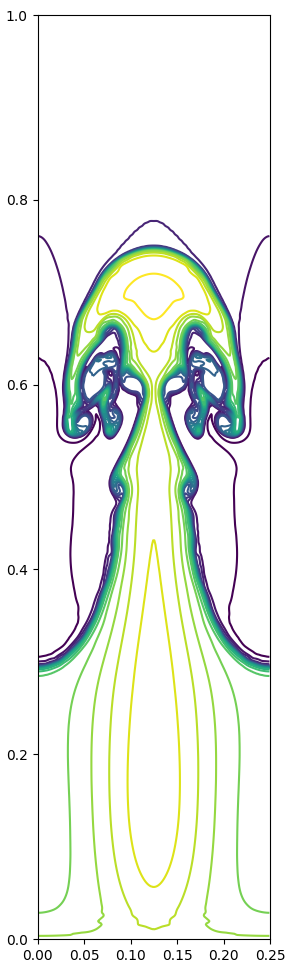}} & \subfloat[WENO5-Z]{\includegraphics[width=0.20\textwidth]{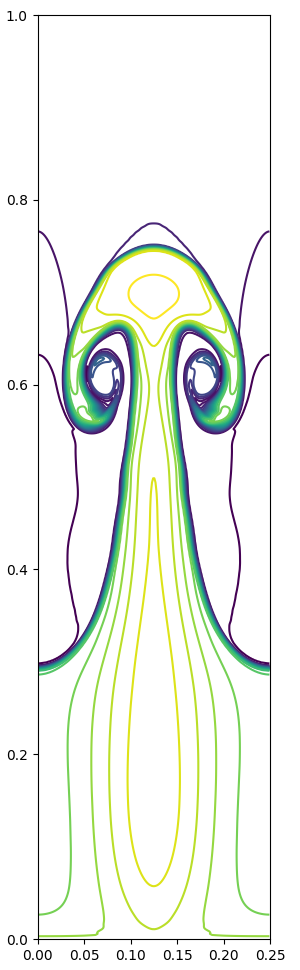}} & \subfloat[WENO5-Z-NN ]{\includegraphics[width=0.20\textwidth]{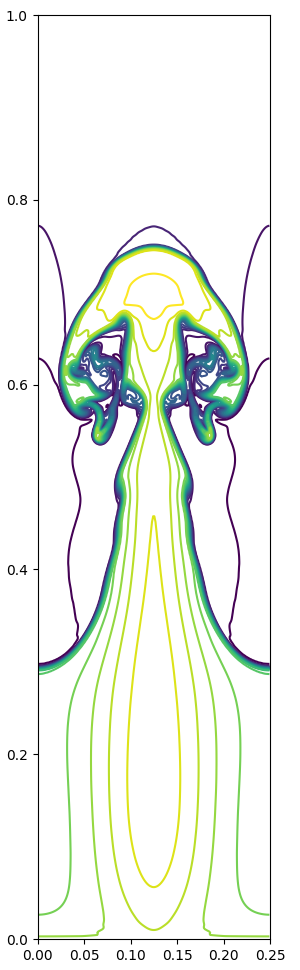}}
		\end{tabular}
		\caption{Density contours of the Rayleigh--Taylor instability problem at $t = 1.95$ with $200 \times 800$ grid points. This figure is drawn with $20$ density contours between $0.9$ and $2.2$.}
		\label{fig:4.4.2}
	\end{figure}

	The Rayleigh--Taylor instability problem \cite{ref49} describes the instability that develops at the interface between two fluids of differing densities when a heavy fluid rushes towards a light fluid with a certain acceleration. It is commonly used to assess the dissipation properties of numerical schemes. The computational domain is $[0, 0.25] \times [0, 1]$ and the specific heat ratio is $\gamma = \frac{5}{3}$. The initial condition is
	\begin{align}
		(\rho, u, v, p) = 
		\begin{cases}
			(2, 0, -0.025\sqrt{\frac{\gamma p}{\rho}}\cos(8\pi x), 2y + 1), & \text{if} \quad y < 0.5,	\\
			(1, 0, -0.025\sqrt{\frac{\gamma p}{\rho}}\cos(8\pi x), y + 1.5), & \text{otherwise}.
		\end{cases}
		\notag	
	\end{align}
	The source term $S = [0, 0, \rho, \rho v]^T$ is added to the right-hand side of the 2D Euler equations in Eq. (\ref{eq:4.4.1}). The reflective boundaries are used on the left and right. The top and bottom boundaries are fixed by $(\rho, u, v, p) = (1, 0, 0, 2.5)$ and $(\rho, u, v, p) = (2, 0, 0, 1)$. Fig. \ref{fig:4.4.2} compares the density contours on the grid $200 \times 800$ at $t=1.95$. 
	Due to the large numerical dissipation, the WENO5-JS and WENO5-Z schemes fail to capture small-scale vortices of the instabilities. In contrast, the WENO5-JS-NN and WENO5-Z-NN schemes exhibit notable improvements, resolving much richer vortical structures.
	Moreover, our schemes retain the symmetry of Rayleigh--Taylor instability during its development.
	
	\noindent \textbf{Example 8.} Double Mach reflection problem.
	
\begin{figure}[htbp]
		\centering
		\centering
		\begin{tabular}{cc}
			\subfloat[WENO5-JS]{\includegraphics[width=0.35\textwidth]{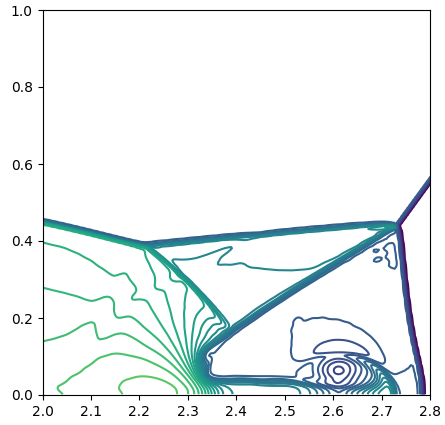}} \hspace{0.1cm} & 	\subfloat[WENO5-Z]{\includegraphics[width=0.35\textwidth]{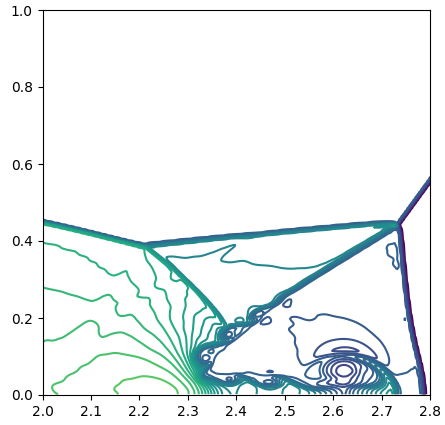}} \\
			\subfloat[WENO5-JS-NN]{\includegraphics[width=0.35\textwidth]{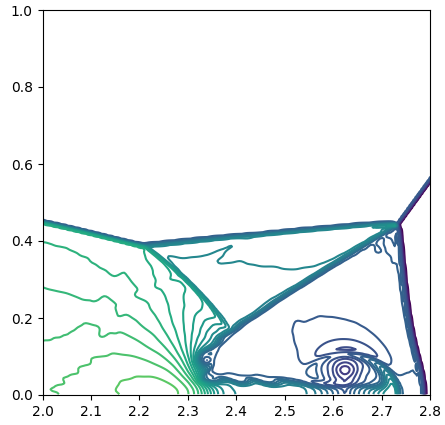}} \hspace{0.1cm} & 	\subfloat[WENO5-Z-NN]{\includegraphics[width=0.35\textwidth]{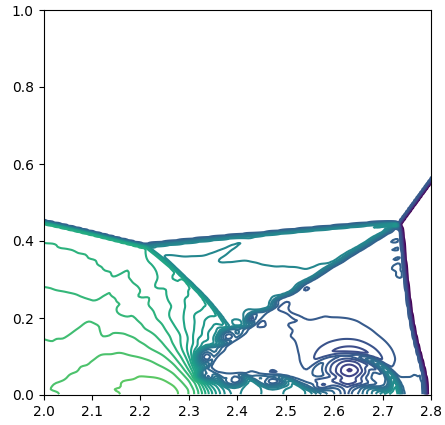}} \\
		\end{tabular}
		\caption{Density contours of the double Mach reflection problem at $t = 0.2$ with $960 \times 240$ grid points. This figure is drawn with $40$ density contours between $2$ and $22$.}
		\label{fig:4.4.3}
	\end{figure}

	This test was first proposed by Woodward and Colella in \cite{ref46} and is used to test the shock-capturing and high-fidelity abilities. The initial condition on the computational domain $[0, 4] \times [0, 1]$ is
	\begin{align}
		(\rho, u, v, p) = 
		\begin{cases}
			(1.4, 0, 0, 1), & x > \frac{1}{6} + \frac{y}{\sqrt{3}},	\\
			(8, 8.25 \sin(60 ^ \circ), -8.25 \cos(60 ^\circ), 116.5), & \text{otherwise}.
		\end{cases}
		\notag	
	\end{align}
	The specific heat ratio is $\gamma=1.4$. The left boundary conditions are post-shocked states, and the right boundary conditions are non-reflective. The top conditions are determined by the exact motion of the oblique shock. The bottom boundary conditions are non-reflective if $x \leq \frac{1}{6}$ and the rest are reflective. The resolution is $960 \times 240$, and the final simulation time is $t = 0.2$. The enlarged contours in the domain $[2, 3] \times [0, 1]$ are presented in Fig. \ref{fig:4.4.3}. The WENO5-Z and WENO5-Z-NN schemes capture rolled-up vortices along the slip line, which are scarcely observed in WENO5-JS and WENO5-JS-NN. The new schemes have higher resolution compared to WENO5-JS/Z.
	
	\section{Conclusion}
	\label{section5}
	In this paper, data-driven optimized high-order WENO schemes are proposed for solving hyperbolic conservation laws. The idea is to introduce a compensation term to the weights of the WENO schemes, which is learned by the neural network. Through the analysis of the ADR, we estimate the spectral error of nonlinear schemes and provide a new approach to improve spectral accuracy. By minimizing the reconstruction error of the WENO5-NN method within a large number of stencils that contain the sine function, the spectral properties of the new schemes are enhanced.
	To inhibit non-physical oscillations near discontinuities and maintain stability at high wavenumbers, the TVD condition and anti-dissipation penalization are applied.
	A series of benchmark cases is used to compare the WENO5-JS/Z-NN with WENO5-JS/Z schemes. The numerical results indicate that the new schemes can capture discontinuous flow structures with higher fidelity and resolve sophisticated multiscale structures with better resolution. The paradigm employed in constructing the WENO5-JS/Z-NN schemes also provides a foundation for developing structure-preserving numerical methods, which we intend to explore in future work.
	
	\section*{Acknowledgements}
		
		Y. Gu was partially supported by the NSFC (No. G0592370101). 
		H. Shen was partially supported by the NSFC (No. 62231016).
		L. Xu was supported by the Innovation Group Grant of Sichuan Province (No. 2025NSFTD0004) and the NSFC (No. 12431015 and No. 62231016). 
		G. Zhou was supported by the NSFC (No. 12171071 and No. 12431015).
		
	\bibliographystyle{elsarticle-num}
	\bibliography{references}

\end{document}